  \documentclass[final,1p,times]{elsarticle}
\usepackage{amssymb}
\usepackage{amsthm}

\usepackage[ruled,vlined,linesnumbered,algoruled,boxed,lined]{algorithm2e}
\usepackage[utf8]{inputenc}
\usepackage[toc,page]{appendix}
\usepackage{amsmath}
\usepackage{graphicx}
\usepackage{grffile}
\usepackage{hyperref,url,color}
\usepackage[linesnumbered,algoruled,boxed,lined]{algorithm2e}
\usepackage[usenames,dvipsnames,svgnames]{xcolor}
\usepackage{tikz}
\usetikzlibrary{arrows,automata,backgrounds,calendar,chains,matrix,mindmap,patterns,petri,shadows,shapes.geometric,shapes.misc,spy,trees,calc,intersections}
\usepackage{wrapfig}
\usepackage{multirow}

\DeclareMathAlphabet\mathbfcal{OMS}{cmsy}{b}{n}

\newcommand{\grad}{\nabla}
\newcommand{\dvg}{\nabla\cdot}

\newcommand{\V}{\left(L^2(\Omega)\right)^{N}}
\newcommand{\A}{\mathcal{A}}
\newcommand{\F}{\mathcal{F}_{\triangle}}
\newcommand{\G}{\mathcal{G}_{\triangle}}
\newcommand{\bfF}{\mathbfcal{F}}
\newcommand{\bfG}{\mathbfcal{G}}
\newcommand{\bfX}{\boldsymbol{\Lambda}}
\newcommand{\x}[1]{{\bf \lambda}_{#1}}
\newcommand{\bfB}{\mathbfcal{B}}
\newcommand{\bfD}{\mathbfcal{D}}
\newcommand{\bfR}{\mathbfcal{R}}

\newcommand{\bfS}{\mathbfcal{S}}

\newcommand{\bfP}{\mathbfcal{P}}

\newcommand{\bfK}{\mathbfcal{K}}

\newcommand{\Hilbert}[2]{\langle#1,#2\rangle_{\V}}
\journal{}

\begin{document}

\begin{frontmatter}



\title{PiTSBiCG: Parallel in Time Stable Bi-Conjugate Gradient Algorithm}

\author{Mohamed Kamel RIAHI}
\address{Department of Applied Mathematics, Khalifa University, PO Box 127788, Abu Dhabi, UAE}
\ead{mohamed.riahi@ku.ac.ae}

\begin{abstract}
This paper presents a new algorithm for the parallel in time (PiT) numerical simulation of time dependent partial/ordinary differential equations. We propose a reliable alternative to the well know parareal in time algorithm, by formulating the parallel in time problem algebraically and solve it using an adapted Bi-Conjugate gradient stabilized method. The proposed Parallel in time Stable Bi-Conjugate algorithm (PiTSBiCG) has a great potential in stabilizing the parallel resolution for a variety of problems. In this work, we describe the mathematical approach to the new algorithm and provide numerical evidences that show its superiority to the standard parareal method.  
\end{abstract}


% \begin{highlights}


\end{frontmatter}


\section{Introduction}
Parallelization across the time directions has been proposed in the early sixsty\cite{Gander50Years}, see also~\cite{ burrage1995parallel} and references therein.  This new fashion of parallelization has been brought up again, to the attention of researchers and engineers with the first appearance of the so-called parareal algorithm by Lions et al~\cite{Maday2001}. Since then, many variants of the of the parareal algorithm have been proposed, tested, and analyzed thoroughly.  Motivated by the flexibility that such algorithm proposes, many applications have attracted the use of the parallel in time method; Among which, Fluid mechanics~\cite{HechtPararealNS}, Neutronics \cite{KamelLWR}, Optimal control~\cite{Maday2013,paraopt}, Quantum physics~\cite{Riahi2016} finance \cite{BalAmericanPut} etc... Interested reader my refer to the recent review~\cite{Gander50Years} and references therein. 

The parareal algorithm benefits from a predictor-corrector scheme to gain acceleration in a shorter time compared to serial computation. In such scheme the prediction is performed by a computationally cheap coarse solver $\G$ while the correction is performed using an accurate high ordered but computationally expensive fine solver $\F$. If the problem is formulated over an interval $[0,T]$ and $(T_n)_{n=0}^{N}$ is a set of times in this interval, then the parareal in time method aims at building a sequence $(\lambda_n^k)_k$ for each time $T_n$ such that, as $k$ goes to infinity, $\lambda_N^k$ converges to the fine solution at time $T_n$. As it will be outlined in section \ref{sec:ParallelSettings}, the sequence $(\lambda_n^k)_k$ is defined through a recursion formula involving predictions of $\G$ and corrections of $\F$. The method decomposes the time domain in the sense that it allows to divide the propagation of the fine solver over $[0,T]$ into propagation on sub-intervals that can be run concurrently and simultaneously on several processors.

\section{Parallel in time settings}\label{sec:ParallelSettings}
The main goal behind time parallelization algorithms is to be able to solve a given time dependent problem over a predefined set of sub-intervals. Ideally, one should be able to solve only on these sub-intervals. Unfortunately, and because the sequential nature of the time-evolution problem, this is not possible unless at least these sub-intervals inter-communicate information. 

For a positive time $T>0$ and a bounded domain $\Omega$ with Lipschitz boundary $\partial\Omega$, we consider the following Cauchy problem 
\begin{equation}\label{Cauchy}\begin{cases}
\dfrac{\partial }{\partial t}y(t,x) -\A(t,y(t,x)) = f(t,x) &\hbox{ in }[0,T]\times\Omega\\
y(0,x)=y_0.&\hbox{ at } \{0\}\times\Omega\end{cases}
\end{equation}
for which $\A(t,\bullet(t,x))$ stands for a differential linear operator, and $f(t,x)$ is a given source term. The Cauchy problem \eqref{Cauchy} is supplemented with appropriate boundary condition at $\partial\Omega$. It is assumed that \eqref{Cauchy} is well posed and has a unique solution on its computational domain. Therefore, it generates a semi-group of evolution that we shall denote by $\F$, in such a way $y(t+\Delta t,x)=\F(y(t,x))$ represents the solution at time $t+\Delta t$ from a given initial condition $y(t,x)$ at time $t$. 
Without loss of generality, we assume that the time interval is equality split into $N$ sub-intervals. The aim of parallel in time computing is solving the problem \eqref{Cauchy} on sub-intervals $[T_{n},T_{n+1}]$ for $n=0,\cdots,N-1$, with $0=T_{0}<T_1<\cdots<T_n<T_{n+1}:=T_{n}+\Delta t<\cdots<T_{N}:=T$. To this end, we attribute $y_n(t,x)$ to the sub-interval $[T_n,T_{n+1}]$, which represents solution to
\begin{equation}\label{Cauchyn}\begin{cases}
\dfrac{\partial }{\partial t}y_{n}(t,x) -\A(t,y_{n}(t,x)) = f_{n}(t,x) &\hbox{ in }[T_{n},T_{n+1}]\times\Omega\\
y_{n}(T_{n},x)=\lambda_{n}&\hbox{ on }\{T_{n}\}\times\Omega.\end{cases}
\end{equation}
For a given initial condition, $\lambda_{n}$, and source term $f_n=f_{|_{[T_{n},T_{n+1}]}}$. Problems \eqref{Cauchyn} are independent and are, hence, solvable in a parallel fashion. Although, compared to the sequential case of \eqref{Cauchy}, the parallel solutions $(y_n)_{n}$ has to satisfy the continuity condition, which writes 
\begin{equation}    \label{continuity}
\lambda_{n+1} = \F\left(y(T_{n},x)\right) = \F\left(\lambda_{n})\right), 
\quad \forall n \in \{0,\cdots,N-1\}
\end{equation}
By writing the collection of the initial conditions $(\x{n})_{n}$ in a vector representation as such $\bfX=(\x{0},\cdots,\x{n},\cdots.\x{N-1})^t$, the continuity conditions is therefore satisfied by the solution to the following algebraic linear system.
\begin{equation}\label{algebraiclinearsystem}
\underbrace{
    \begin{bmatrix}
    I & 0 & \cdots & \cdots & 0\\
    -\F & I &0 & \cdots & \vdots\\
    0 & \ddots & \ddots & \ddots & \vdots\\
    \vdots & \ddots &-\F & I & 0\\
    0  & \cdots &0 &-\F & I 
    \end{bmatrix}}_{{\bfF}}
    \underbrace{\begin{bmatrix}
    \lambda_{0}\\\vspace{0.1in}\vdots\\\lambda_{n}\\\vdots\\\lambda_{N-1}\end{bmatrix}}_{\bfX}
    =\underbrace{\begin{bmatrix}
    \lambda_{0}\\\vspace{0.05in}\vdots \\ \vdots \\\vdots \\0
    \end{bmatrix}.}_{\bfB}  
\end{equation}

In~\cite{Maday2001} the parareal algorithm updates the initial conditions using a predictor-corrector scheme as follows 
\begin{equation}\label{pararealscheme}
    \lambda_{n+1}^{k+1} = \G(\lambda_{n}^{k+1}) +\F(\lambda_{n}^{k+1}) -\G(\lambda_{n}^{k+1})
\end{equation}
where $\F$ and $\G$ are respectively a fine an coarse solver for the time-evolution problem at hand. It has been shown at the early stage of the parareal algorithm that the scheme \eqref{pararealscheme} is algebraically equivalent to 
\begin{equation}\label{pararealAlgebraic}
    \bfX^{k+1} = \bfX^{k} + \bfG^{-1}\left( \bfB - \bfF\cdot\bfX^{k} \right)
\end{equation}
In the present work, we shall investigate the algebraic structure of the non-symmetric operator-block linear system arising from the time domain decomposition. We, indeed, propose new and robust time parallel algorithm that outperform the parareal algorithm. The new method has also the capability of being applied to any evolution equation, as it is the case for the parareal method. On the other hand, our algorithm performs well with large sub-domains compared to the parareal method.   

\section{PiTSBiCG method}
The approach we follow in developing our numerical algorithm is based on the nature of the block-non-symmetric linear system \eqref{algebraiclinearsystem}. It is well-known in the computational linear algebra literature (see for instance~\cite{saad_2003}) that the most appropriate and efficient linear solvers are BiCG/QMR and their enhanced variant such as GMRes. Our focus goes for the BiCG which uses the Lanczos Biorthogonalization procedure, and process the solution through projection upon the following Krylov vector space 
$${\bf\mathcal{K}}_{m}=\left\{v_{1},\bfF\cdot v_{1},\cdots,\bfF^{m-1}\cdot v_{1}\right\},$$
and orthogonal to 
$${\bf\mathcal{L}}_{m}=\left\{w_{1},(\bfF^T)\cdot w_{1},\cdots,(\bfF^T)^{m-1}\cdot w_{1}\right\},$$
It is worth mentioning that, minimal residual algorithms, such as QMR and GMRES, that are based on matrix algebra factorization, are not suitable to solve \eqref{algebraiclinearsystem} as of the operators-block-structure of the problem. In the solution process of \eqref{algebraiclinearsystem} we only consider the matrix-by-vector product., where the operator block structures are neither assembled nor stored. Furthermore, algorithms that do not use factorization techniques can, actually, be used to solve the parallel in time problem \eqref{algebraiclinearsystem}. Although, this may potentially engage operator-block-structure matrix transpose calculation in a least square CG-based programming~\cite{RAO20121021,chen2015adjoint}, or with the use of the BiCG method. Nonetheless, the CGS algorithm which is designed to avoid evaluating the transpose of the linear system my also face accuracy challenges as such in the matrix linear algebra, where the residual may present unacceptable high variations hence affects the outcome in the iterative process, where substantial buildup of rounding error may often be observed. The BiCGStab combines residual formulas from both BiCG and CGS and has been shown to be effective in both classical and block versions~\cite{Simoncini1997}. The later observation encouraged us to consider such algorithm for the time parallel settings.

In the sequel, we adapt the well known linear algebra iterative solver BiCGStab to the time parallelisation setting and consider solving the ``fictive\footnote{fictive: means that the linear system is neither assembled nor stored.}'' linear system $\bfF\cdot\bfX=\bfB$, which we precondition using the left-preconditioner operator-block matrix $\bfG^{-1} \approx \bfF^{-1}$. Indeed, the inverse of the matrix $\bfF$ writes simply 
\begin{equation}\label{eqinverse}
\bfF^{-1}:=
    \begin{bmatrix}
    I & 0 & \cdots & \cdots & 0\\
    \F & I &0 & \cdots & \vdots\\
    \F^{2} & \ddots & \ddots & \ddots & \vdots\\
    \vdots & \ddots &\F^{} & I & 0\\
    \F^{N}  & \cdots &\F^{2} &\F & I 
    \end{bmatrix},\,
    \begin{array}{c}
         \text{ hence }\\
         \text{approximated}\\
         \text{by}
    \end{array}, \quad    \bfG^{-1}:=
    \begin{bmatrix}
    I & 0 & \cdots & \cdots & 0\\
    \G & I &0 & \cdots & \vdots\\
    \G^{2} & \ddots & \ddots & \ddots & \vdots\\
    \vdots & \ddots &\G^{} & I & 0\\
    \G^{N}  & \cdots &\G^{2} &\G & I 
    \end{bmatrix}
\end{equation}

Therefore, in practice we solve 
\begin{equation}\label{preconditionedlinearsystem}
    \bfG^{-1}\cdot\bfF\cdot\bfX=\bfG^{-1}\cdot\bfB
\end{equation}
without of course assembling this linear system. We rather use the matrix-by-vector output products, following the steps of {\bf Algorithm}~\ref{algo}.   

\begin{wrapfigure}{l}{0.6\textwidth}
\begin{algorithm}[H]
\DontPrintSemicolon
\SetAlgoLined
 \KwIn{Ininital guess $\bfX^{0}=(\x{0},\cdots,\x{N-1})^T$,\\ tolerance $\epsilon$, restart tolerance $\epsilon_{0}$\;}
 $\bfR^{0}\gets\bfG^{-1}(\bfB-\bfF\cdot \bfX^{0})$\;
 $\tilde\bfR\gets\bfR^{0}$\;
 $\bfP^{0}\gets\bfR^{0}$\;
 $k\gets0$\;
 \While{$\|\bfR^{k}\|_{N,\infty}>\epsilon$}
 {
 \vbox{\colorbox{white}{$\bfD^{k}\gets\bfG^{-1}\cdot\bfF\cdot\bfP^{k}$\;}}
 $\alpha_{k} = \Hilbert{\bfR^{k}}{\tilde\bfR}\slash \Hilbert{\bfD}{\tilde\bfR}$\;
 $\bfS^{k}\gets \bfR^{k} - \alpha_{k}\bfD^{k}$\;
 \vbox{\colorbox{white}{$\bfK^{k}\gets\bfG^{-1}\cdot\bfF\cdot\bfS^{k}$\;}}
 $\omega_{k}=\Hilbert{\bfK^{k}}{\bfS^{k}} \slash \Hilbert{\bfK^{k}}{\bfK^{k}}$\;
 $\bfX^{k+1}\gets \bfX^{k} +\alpha_{k}\bfP^{k}+\omega_{k}\bfS^{k}$\;
 $\bfR^{k+1} \gets \bfS^{k} - \omega_{k}\bfK^{k}$\; 
 $\beta_{k} = \frac{\alpha_{k}}{\omega_{k}} \Hilbert{\bfR^{k+1}}{\tilde\bfR}\slash \Hilbert{\bfR^{k}}{\tilde\bfR}$\;
 $\bfP^{k+1} \gets \bfR^{k+1} + \beta_{k} \left(\bfP^{k} - \omega_{k}\bfD^{k} \right)$\;
      \If{$\left|\Hilbert{\bfR^{k+1}}{\tilde\bfR}\right|\leq \epsilon_{0}$}
    {$\tilde\bfR\gets\bfR^{k+1}$\;
    $\bfP^{k}\gets\bfR^{k+1}$\;}
 $k\gets k+1$\;
 }
 \KwResult{ $\bfX^{\star}=\bfX^{k+1}$\;}
 \caption{\hbox{PiTSBiCG algorithm}}\label{algo}
\end{algorithm}
\end{wrapfigure}

The pseudo-code is depicted in {\bf Algorithm}~\ref{algo}, in which we consider the PiTSBiCG to solves the parallel across the time direction linear system in its preconditioned version \eqref{preconditionedlinearsystem}. The structure of the algorithm is very much close to the linear algebraic version of the BiCGStab well-known algorithm, although we emphasize the block-operator structure issued from the parallel across the time direction problem.  Assume we have $N$ sub-domain for the parallel in time computation, the PiTSBiCG generates $L^{2}(\Omega)$ solutions through its iterative process that are defined at the breakpoints. Algebraic-wise, the set of all solutions produced, whether they are initial conditions, residuals, or directions they are piledup into a vector form element of the vector space $\left(L^{2}(\Omega)\right)^{N+1}$. Such vector space is endowed with the following block-structured inner product $\langle\cdot,\cdot\rangle_{\left(L^{2}(\Omega)\right)^{N}}=\sum_{n=0}^{N-1}
\langle\cdot,\cdot\rangle_{L^{2}(\Omega)}$  . 
It is worth noting that the inverse symbol applied for the Matrix-operator $\bfG$ is described as a preconditioner at lines 6 and 9, which means that we solve the preconditioning resolution sequentially.  This is, permitted as the semi-group operators $\G$ are made, by construction, coarse and non-expansive. One can see from \eqref{eqinverse} that such inverse can easily be solved following forward substitution and this is nothing but a sequential resolution. We here recall that none of the operator-matrices is stored, and only a matrix-by-vector product is involved in our algorithm. 
We present in the sequel applications of our method in a variety of examples and compare its performance with the one of the parareal algorithm.

\section{Numerical evidences}

We consider the following advection-diffusion-reaction equation
\begin{equation}\label{eqtest}
\begin{cases}
     \dfrac{\partial }{\partial t}y(t,x) + {\bf u}\cdot\nabla y(t,x) + \dvg(\mu(t,x)\grad y) + r(t,x) y(t,x) = f(t,x) &\hbox{ in }[0,T]\times\Omega\\
     y(t,x) = g(t,x) &\hbox{ in }[0,T]\times\partial\Omega\\
     \dfrac{\partial y}{\partial \vec{n}}(t,x) = h(t,x) &\hbox{ in }[0,T]\times\partial\Omega\\
     y(0,x) = y_0 &\hbox{ at }\{0\}\times\Omega\\
\end{cases}
 \end{equation}
\begin{itemize}
     \item Diffusion equation: in \eqref{eqtest} we consider ${\bf u}=0$, $r(t,x)=0$, and $\mu=1$.\vspace{-.1in}
     \item Diffusion reaction equation: in \eqref{eqtest} we consider ${\bf u}=0$, $r(t,x)=1.5$, and $\mu=1$.\vspace{-.1in}
     \item Advection diffusion reaction: in \eqref{eqtest} we consider ${\bf u}=[-y,x]$. We put $r(t,x)=5\cdot10^{-1}$ and reduce the diffusion term $\mu=1\cdot10^{-1}$.
 \end{itemize}

\begin{figure}[htbp]\centering
\begin{tabular}{ccc}
&Iterations & Matrix-vector multiplication\\
\rotatebox{90}{\hspace{1in}Diffusion}&
\includegraphics[width=5cm,height=5cm]{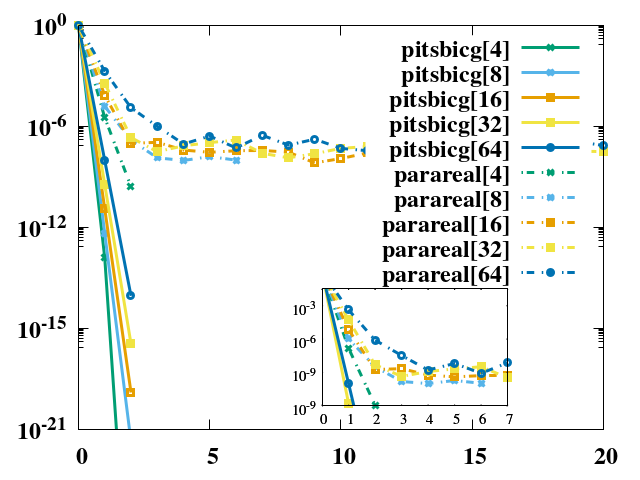} &
\includegraphics[width=5cm,height=5cm]{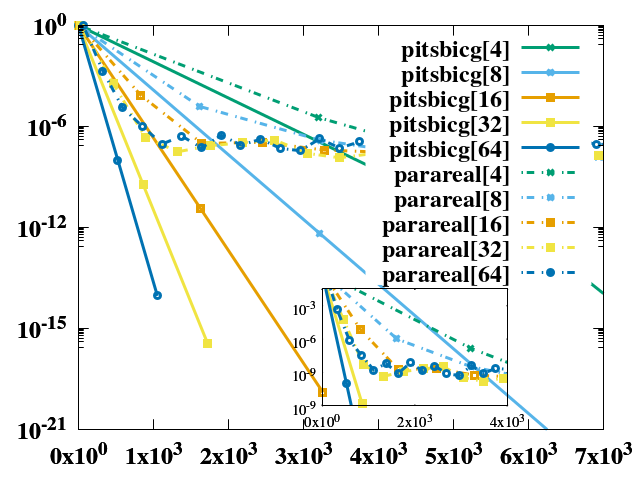}\\   
\rotatebox{90}{\hspace{0.6in}Diffusion Reaction}&
\includegraphics[width=5cm,height=5cm]{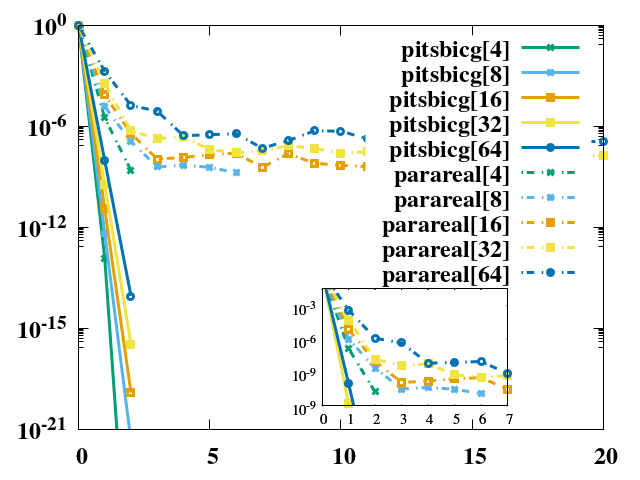} &
\includegraphics[width=5cm,height=5cm]{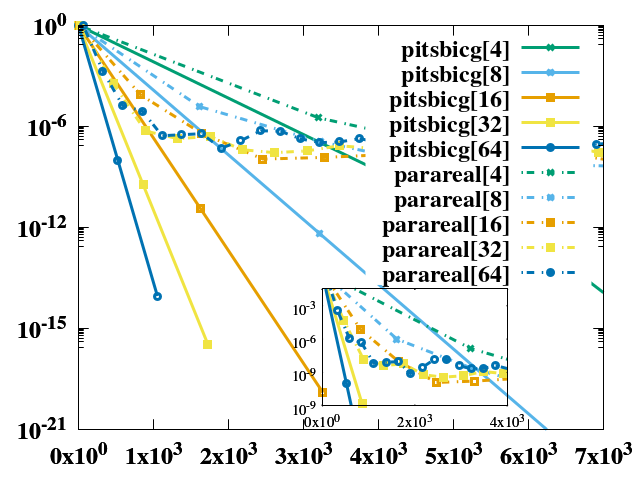}\\
\rotatebox{90}{\hspace{.3in}Advection Reaction Diffusion}&
\includegraphics[width=5cm,height=5cm]{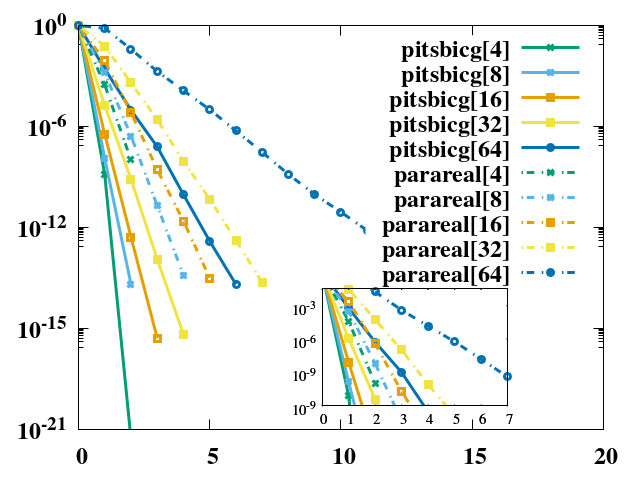} &
\includegraphics[width=5cm,height=5cm]{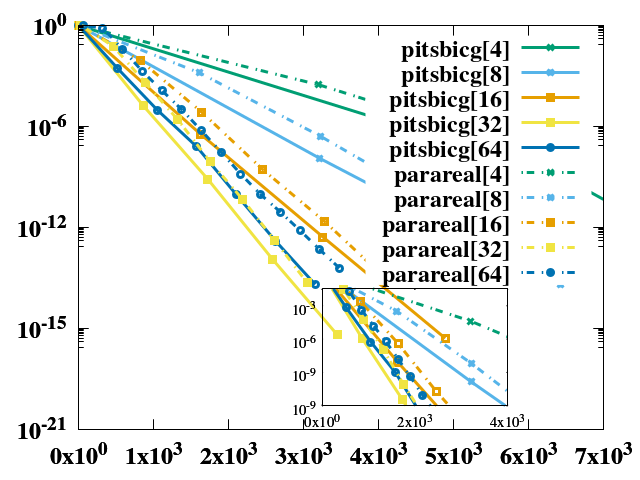}\\
\end{tabular}
    \caption{Residual decreasing in values versus number of iterations (left column) and versus Matrix-by-vector multiplication (right column), Diffusion PDE case (first row), Diffusion reaction PDE case (second row), and Advection reaction diffusion PDE case (third row).}
    \label{fig:results}
\end{figure}

The results of our numerical simulations are depicted in Figure~\eqref{fig:results}. We used the Finite Element method for the space discretize the handled PDEs with step size $h=2\cdot10^{-2}$, where the computational domain $\Omega$ is taken to be a unit square centered at the origin. We also used Backward Euler time marching for the time discretization with step size $dt=1.10^{-3}$, where the total physical time $T=6.4$. The initial condition $y_0$ is taken as a Gaussian.  These settings applies for all numerical simulation we have conducted, using FreeFem++~\cite{freefem} software.  

For each type of the PDEs in \eqref{eqtest}, we consider $2^n$ sub-domains with $n=2,3,\cdots,6$. These sub-domains are allocated to $2^n$ CPUs in a Linux workstation. The results are presented in term of the decreasing in value of the residual of the preconditioned linear system versus the iterations (left) and versus the matrix-by-vector multiplication (right). The results show the clear superiority of the proposed algorithm compared to the plain version of the parareal algorithm. The clear out-performance is shown in both the number of iterations and in terms of the operations of matrix-by-vector multiplication. We note here that PiTSBiCG makes two iterations (Bi-directions) compared to the parareal algorithm, this means if one wants to evaluate the iterations that ensure decent in value of the residual has to multiply by two the iteration of PiTSBiCG. For this reason, it is more appropriate to analyze the performance in terms of the Matrix-by-vector product.

Future consideration, as extension of the presented results, includes the analysis of the algorithm and exploit its applicability in real-world engineering's problems involving differential equation whether ordinary, partial or fractional. 
\section{Conclusion}

This paper presented a novel parallel in time algorithm based on simple investigation of the equivalent algebraic structure of the initial condition system of equations. We showed that our method outperforms the well-known parareal method and presented several numerical examples to support our claim.



  \bibliographystyle{elsarticle-num} 
  \bibliography{biblio.bib}


% \end{thebibliography}
\end{document}